\documentclass[12pt]{article}
\usepackage{amsmath, amssymb, amsthm, bm}
\usepackage{graphicx}
\usepackage{authblk}

\title{\bf Grothendieck's Geometric Universes and A Sheaf-Theoretic Foundation of Information Network}

\author{\Large Takao Inou\'{e}}

\affil{\large Faculty of Informatics, Yamato University, \\ Osaka, Japan\footnote{Email: inoue.takao@yamato-u.ac.jp; \\ Personal Email: takaoapple@gmail.com \\ [I prefer my personal email address for correspondence.]}} 

\date{February 19, 2026}

\newtheorem{definition}{Definition}[section]

\begin{document}
\maketitle

\begin{abstract}
This paper proposes an interpretation of Grothendieck's geometric universes as a foundational framework for \emph{information networks}. 
We argue that Grothendieck topologies, sheaves, and topoi provide a sheaf-theoretic semantics in which distributed and locally held information can be integrated into globally coherent structures. 
In this setting, local informational states are represented by sections, while the sheaf condition governs consistency, agreement, and consensus across a network. 
Logical validity and mathematical existence are therefore not imposed externally but arise intrinsically from geometric and categorical conditions. 
From this perspective, Grothendieck's geometric universes constitute a natural foundation for information networks governed by intrinsic logical principles. 
Moreover, we propose that Grothendieck's geometric universes themselves concretely instantiate what the author calls \emph{intrinsic logicism}. 
This position is intended as a contemporary reconstruction of the classical logicist program of Frege and Russell, reformulated within the framework of category theory and topos theory, where logical structure is generated internally by geometric and categorical organization rather than presupposed as an external foundational layer.
\\ \medskip

\\
\medskip

\end{abstract}

\tableofcontents

\section{Introduction}

Grothendieck's revolution in mathematics replaced point-based geometric intuition with a structural and relational perspective \cite{GrothendieckSGA4,JohnstoneTopos}. 
Beyond its impact on geometry, this shift provides a powerful conceptual framework for understanding \emph{information networks}, in which information is distributed, contextual, and locally constrained. 
Grothendieck topologies, sheaves, and topoi offer a formal language for describing how such locally available information can be communicated, compared, and coherently integrated.

In this paper, we argue that Grothendieck's geometric universes should be understood as a sheaf-theoretic foundation of information networks. 
Local sections correspond to partial informational states held by individual nodes, while the sheaf condition formalizes principles of consistency, agreement, and consensus across a network. 
From this viewpoint, logical validity is not an external constraint imposed on information processing, but an intrinsic consequence of the geometric and categorical structures governing information flow, in the spirit of Lawvere's categorical foundations \cite{LawvereLogic}.

Moreover, we propose that this intrinsic emergence of logical validity reflects a foundational position that the author calls \emph{intrinsic logicism}. 
This position may be regarded as a contemporary reconstruction of the classical logicist program of Frege and Russell: rather than grounding mathematics in an externally fixed logical system, logical structure is understood as arising internally from the geometric and categorical organization of mathematical universes themselves. 
In this sense, Grothendieck's geometric universes do not merely support logical reasoning but concretely instantiate a modern form of logicism adapted to distributed, networked, and structurally mediated contexts.

\section{Grothendieck Topologies and Sheaves}
We recall the precise definition of a Grothendieck topology, which is central to the subsequent development and therefore retained in full detail.

\begin{definition}
Let $\mathcal{C}$ be a category. A \emph{Grothendieck topology} $J$ on $\mathcal{C}$ assigns to each object $U$ of $\mathcal{C}$ a collection $J(U)$ of families of morphisms $\{U_i \to U\}_{i \in I}$, called \emph{covering families}, satisfying the following axioms:
\begin{enumerate}
\item (Isomorphism) If $V \xrightarrow{\cong} U$ is an isomorphism, then ${V \to U} \in J(U)$.
\item (Stability under pullback) If ${U_i \to U} \in J(U)$ and $V \to U$ is any morphism, then the pullback family ${U_i \times_U V \to V}$ belongs to $J(V)$.
\item (Transitivity) If $\{U_i \to U\}_{i \in I} \in J(U)$ and for each $i \in I$,
$\{V_{ij} \to U_i\}_{j \in J_i} \in J(U_i)$, then the composed family
$\{V_{ij} \to U\}_{i \in I,\, j \in J_i}$ belongs to $J(U)$.

\end{enumerate}
\end{definition}

A \emph{sheaf} on the site $(\mathcal{C},J)$ is a presheaf $F : \mathcal{C}^{op} \to \mathbf{Set}$ satisfying locality and gluing conditions with respect to the covering families of $J$. This definition abstracts local consistency and global coherence beyond topological spaces.

The author believes that the definition of Grothendieck topology itself has an affinity with the structure of information networks.

\section{Topoi as Geometric Universes}
A Grothendieck topos may be understood as a generalized geometric universe in which both mathematical objects and logical reasoning are internalized \cite{MacLaneMoerdijk,JohnstoneSketches}. A topos possesses finite limits, exponentials, and a subobject classifier, thereby supporting an internal higher-order intuitionistic logic. From the perspective of intrinsic logicism, a topos is not merely a model of logic but a structure in which logic is generated intrinsically by geometric conditions.

\section{Intrinsic Logicism and Internal Semantics}

The intrinsic logic of a topos is formulated via its internal language, in which propositions correspond to subobjects and logical operations are interpreted categorically \cite{MacLaneMoerdijk,LawvereLogic}. 
In this setting, logical reasoning is carried out internally to a given geometric universe, and truth is evaluated relative to its structural and contextual conditions. 
Logical validity is therefore local, contextual, and structurally constrained, rather than absolute in the sense of an externally fixed logical calculus.

This internal perspective motivates what the author calls \emph{intrinsic logicism}. 
Intrinsic logicism is the foundational thesis that logical structure is not imposed on a mathematical universe from an external set-theoretic or syntactic framework, but instead emerges from the internal geometric and categorical organization of that universe itself. 
In contrast to classical logicism, which sought to reduce mathematics to a predetermined logical system, intrinsic logicism reinterprets the logicist aspiration in a structural and geometric form.

From this viewpoint, Grothendieck topoi provide concrete realizations of intrinsic logicism. 
Their internal logics arise naturally from the sheaf-theoretic and categorical conditions that govern how local data are assembled into global structures. 
In particular, the interpretation of propositions as subobjects and of logical operations as categorical constructions reflects the manner in which informational consistency and coherence emerge within a networked system.

Accordingly, intrinsic logicism may be regarded as a contemporary reconstruction of the Frege--Russell logicist program, reformulated in terms of internal semantics rather than external reduction. 
Logical principles are preserved not as universal axioms imposed from outside, but as invariant features generated within geometric universes themselves. 
This shift renders logicism compatible with distributed, contextual, and dynamically structured domains such as information networks.

\section{Sheaf Semantics and Information Networks}
Sheaf semantics provides a natural mathematical model of distributed and contextual information \cite{MacLaneMoerdijk}. From the standpoint of intrinsic logicism, Grothendieck's topos theory may be reinterpreted not only as a geometric framework but also as a logical protocol for managing \emph{information coherence} and \emph{consensus} in distributed systems.

\subsection*{Sheaf Axioms as an Integration Algorithm for Distributed Information}
The axioms of a sheaf, in particular the \emph{gluing condition}, specify how locally held fragments of data across a network can be integrated into a globally consistent form of truth.

\begin{itemize}
\item \textbf{Local data (local sections).} For each node $U$ in the network, a section $s \in F(U)$ represents information available only within the observational or contextual scope of that node.
\item \textbf{Restriction maps as communication.} A restriction morphism $\mathrm{res}_{U,V}$ corresponds to the transmission of information through a communication channel, translating data into a shared context over an overlap $U \cap V$ so that it can be meaningfully compared or combined with information from neighboring nodes.
\item \textbf{Guarantee of consensus.} The \emph{existence} and \emph{uniqueness} clauses of the sheaf axioms ensure that locally compatible pieces of information, when integrated, determine a single global section without redundancy or contradiction. This formalizes consensus formation within the network.
\end{itemize}

In this way, sheaf semantics provides a rigorous foundation for information networks in which global meaning is not imposed externally but emerges intrinsically from the coherent interaction of local informational states. Such networks exemplify intrinsic logicism in practice: logical consistency arises from internal structural conditions rather than from an overarching external truth predicate.

Sheaf semantics provides a natural mathematical model of distributed and contextual information \cite{MacLaneMoerdijk}. Local sections represent partial information, while the sheaf condition enforces global coherence. From the viewpoint of intrinsic logicism, this semantic structure exemplifies how logical consistency is generated internally through gluing conditions rather than imposed by a global truth predicate. Such a perspective is particularly well-suited to information networks and contextual reasoning systems.

\section{Conclusion and Future Work}

In this paper, we have articulated and defended \emph{intrinsic logicism} as a foundational standpoint according to which logical structure is generated internally by geometric and categorical organization. 
Through Grothendieck topologies, sheaves, and topoi, mathematical universes acquire an intrinsic logical semantics that reflects locality, coherence, and contextual validity. 
From this perspective, sheaf theory provides not merely a technical tool but a conceptual foundation in which logical validity and mathematical existence arise from internal structural conditions.

This sheaf-theoretic viewpoint clarifies the internal logic of modern geometry and, at the same time, offers a natural semantic framework for \emph{information networks}. 
Distributed informational states, local consistency conditions, and global coherence can be understood uniformly within the internal semantics of geometric universes. 
In this sense, intrinsic logicism provides a principled account of logic suited to distributed, networked, and context-sensitive systems.

As future work, we aim to develop concrete applications of this paradigm to real-world information networks. 
In particular, we plan to investigate how sheaf-theoretic semantics and intrinsic logical principles can be applied to knowledge representation, consistency management, and consensus formation in distributed informational environments.

$$ $$

\noindent Takao Inou\'{e}

\noindent Faculty of Informatics

\noindent Yamato University

\noindent Katayama-cho 2-5-1, Suita, Osaka, 564-0082, Japan

\noindent inoue.takao@yamato-u.ac.jp
 
\noindent (Personal) takaoapple@gmail.com (I prefer my personal mail)

\bigskip

\end{document}